\documentclass[a4paper,12pt]{article}

\usepackage{amsfonts}
\usepackage{amssymb}
\usepackage{amsbsy}
\usepackage{amsmath}

\begin{document}

\title{Solubility Criteria for Hopf-Galois Structures}

\author{Nigel P.~Byott, \\ 
Department of Mathematics, \\
College of Engineering, Mathematics and Physical Sciences, \\
University of Exeter, 
Exeter EX4 4QF, 
UK\\
E-mail: N.P.Byott@ex.ac.uk\\
}

\newcommand{\Z}{{\mathbb{Z}}}
\newcommand{\Q}{{\mathbb{Q}}}
\newcommand{\F}{{\mathbb{F}}}
\newcommand{\PP}{{\mathbb{P}}}

\newcommand{\gp}{{\mathfrak{p}}}
\newcommand{\gP}{{\mathfrak{P}}}

\newcommand{\sP}{{\mathcal{P}}}
\newcommand{\sW}{{\mathcal{W}}}

\newcommand{\Gal}{\mathrm{Gal}}
\newcommand{\Hom}{\mathrm{Hom}}
\newcommand{\Hol}{\mathrm{Hol}}
\newcommand{\Perm}{\mathrm{Perm}}
\newcommand{\id}{\mathrm{id}}
\newcommand{\lra}{\longrightarrow}
\newcommand{\bs}{\backslash}

\newtheorem{theorem}{\sc Theorem}
\newtheorem{proposition}[subsection]{\sc Proposition}
\newtheorem{corollary}{\sc Corollary}
\newtheorem{lemma}[subsection]{\sc Lemma}

\newenvironment{remark}{\begin{rk}\rm}{\end{rk}} 
\newtheorem{rk}[subsection]{\sc Remark}

\newenvironment{example}{\begin{ex}\rm}{\end{ex}} 
\newtheorem{ex}[subsection]{\sc Example}

\newenvironment{defn}{\begin{df}\rm}{\end{df}} 
\newtheorem{df}[subsection]{\sc Definition}

\newcommand{\qed}{\hspace*{\fill} \rule{2mm}{3mm}}

\newenvironment{pf}{{\em Proof. }}{\qed \newline}
\newenvironment{pf-of}[1]{{\noindent\em Proof of #1.}}{\qed \newline}

\newcommand{\Aut}{\mathrm{Aut}}
\newcommand{\Inn}{\mathrm{Inn}}
\newcommand{\Out}{\mathrm{Out}}
\newcommand{\GL}{\mathrm{GL}}
\newcommand{\SL}{\mathrm{SL}}
\newcommand{\PSL}{\mathrm{PSL}}
\newcommand{\PGL}{\mathrm{PGL}}
\newcommand{\PO}{\mathrm{P}\Omega}
\newcommand{\PSp}{\mathrm{PSp}}
\newcommand{\PSU}{\mathrm{PSU}}
\newcommand{\SU}{\mathrm{SU}}

\newcommand{\eab}{e_\mathrm{ab}}
\newcommand{\enil}{e_\mathrm{nil}}

\newcommand{\N}{\mathbb{N}}

\newcommand{\onto}{\twoheadrightarrow}

\maketitle{}

\centerline{2010 Mathematics Subject Classification: 12F10, 16T05}

\begin{abstract}
Let $L/K$ be a finite Galois extension of fields with group $\Gamma$.
Associated to each Hopf-Galois structure on $L/K$ is a group $G$ of
the same order as the Galois group $\Gamma$. The {\em type} of the
Hopf-Galois structure is by definition the isomorphism type of $G$.
We investigate the extent to which general properties of either of the
groups $\Gamma$ and $G$ constrain those of the other. Specifically, we
show that if $G$ is nilpotent then $\Gamma$ is soluble, and that if
$\Gamma$ is abelian then $G$ is soluble. The proof of the latter
result depends on the classification of finite simple groups. In
contrast to these results, we give some examples where the groups
$\Gamma$ and $G$ have different composition factors. In particular, we
show that a soluble extension may admit a Hopf-Galois structure of
insoluble type.
\end{abstract}

\section{Introduction and statement of results}

Hopf-Galois theory was initiated by Chase and Sweedler \cite{CS},
motivated in part by a wish to develop a version of Galois theory for
inseparable field extensions. Their approach nevertheless gives an
interesting perspective on the classical theory for separable
extensions. The problem of finding all Hopf-Galois 
structures on a given separable extension $L/K$ was expressed in
terms of group theory by Greither and Pareigis \cite{GP}, who showed that 
$L$ may admit $H$-Galois structures for a number of different $K$-Hopf
algebras $H$. Moreover, it can happen that the same Hopf algebra $H$ may have several different actions on $L$, giving rise to different
Hopf-Galois structures \cite{CCo, CRV-JA}.  

In the special case
that $L/K$ is a Galois extension (i.e.~normal as well as separable) with group $\Gamma=\Gal(L/K)$,
the main result of \cite{GP} is that the Hopf-Galois structures on $L/K$ correspond to
regular subgroups $G$ of the group $\Perm(\Gamma)$ of permutations of
$\Gamma$ such that $G$ is normalized by left translations by
$\Gamma$. The
groups $\Gamma$ and $G$ necessarily have the same order, but in
general need not be isomorphic.  We refer to the isomorphism type of
$G$ as the {\em type} of the Hopf-Galois structure.

A number of authors have used the framework developed in \cite{GP} to investigate Hopf-Galois 
structures on various classes of separable extensions.  Apart from their intrinsic interest, one motivation for studying multiple 
Hopf-Galois structures on the same extension is their relevance to questions of integral
Galois module structure: the ring of integers in a Galois
extension of local fields may have better module-theoretic properties
in one of the nonclassical Hopf-Galois structures than it does in the
classical Galois structure \cite{assoc}.

A considerable amount is now known about Hopf-Galois structures on
Galois extensions $L/K$. For an odd prime $p$, a cyclic extension of
degree $p^n$ admits precisely $p^{n-1}$ Hopf-Galois structures, all of
cyclic type \cite{Kohl}, while an elementary abelian extension of
degree $p^n$ admits at least $p^{n(n-1)-1}(p-1)$ Hopf-Galois
structures of elementary abelian type if $p>n$ \cite{Ch05}.  If $p$
and $q$ are distinct primes such that there exists a nonabelian group
of order $pq$, then any Galois extension of degree $pq$, whether
abelian or not, admits Hopf-Galois structures both of abelian type and
of nonabelian type \cite{pq}. More generally, Galois extensions of
degree $mp$, with $p$ prime and $m<p$, are considered in
\cite{kohl-mp}.  In contrast to the cyclic prime-power case, ``most''
abelian extensions admit Hopf-Galois structures of nonabelian type
\cite{BC}. A Galois extension whose group $\Gamma$ is a nonabelian
simple group admits precisely two Hopf-Galois structures, both of type
$\Gamma$ \cite{simple}, whereas an extension with symmetric Galois
group $S_n$ admits many Hopf-Galois structures of type $S_n$ and also
many of type $A_n \times C_2$ \cite{CaC}.

There are fewer results relating to the more general situation where
$L/K$ is separable but not necessarily normal. A separable extension
of prime degree is Hopf-Galois if and only if its Galois closure has
soluble Galois group \cite{Ch-p}, and, for an odd prime $p$, the
Hopf-Galois structures on a radical extension $K(\sqrt[p^n]{a})/K$ of
degree $p^n$ are enumerated in \cite{Kohl}, the result depending on
which roots of unity occur in $K$. The Hopf-Galois characters of all
separable extensions $L/K$ of degree $\leq 6$ have recently been
determined \cite{CRV-lattice}, together with those of the extensions
$F/K$ with $L \subset F \subset N$, where $N$ is the normal closure of
$L/K$.

For a more extensive review of results in separable Hopf-Galois
theory, we refer the reader to \cite{CRV-survey}.

In this paper, we consider only Galois extensions. Let $L/K$ be a
Galois extension with Galois group $\Gamma$, and suppose that $L/K$
admits a Hopf-Galois structure of type $G$. We investigate the extent
to which group-theoretic properties of either one of the groups
$\Gamma$ or $G$ constrain the other. Already in the case of groups of
order $pq$, mentioned above, it may happen that either one of the
groups is abelian while the other is not (and indeed is not
nilpotent). This shows, for example, that it is possible for an
extension of finite fields to admit a Hopf-Galois structure which is
not of nilpotent type. Our two main results will give criteria, in terms of each
of the groups, for the other to be soluble:

\begin{theorem} \label{nilp-sol}
With the above notation, if $G$ is nilpotent then $\Gamma$ is soluble. Thus,
if a finite Galois extension of fields admits a Hopf-Galois structure of
nilpotent type, then the extension has soluble Galois group.
\end{theorem}

\begin{theorem} \label{ab-sol}
If $\Gamma$ is abelian then $G$ is soluble. Thus, any Hopf-Galois
structure on a finite abelian field extension has soluble type.
\end{theorem}

Theorem \ref{nilp-sol} is an application of the theory
of Hall $p'$-subgroups. The proof of Theorem \ref{ab-sol} depends on the
classification of finite simple groups, and in particular uses a
result of Vdovin \cite{V} which bounds the size of abelian
subgroups in a nonabelian simple group.    

In the situations of Theorems \ref{nilp-sol} and \ref{ab-sol},
$\Gamma$ and $G$ are soluble groups of the same order, so certainly
have the same composition factors. In the various results for Galois
extensions $L/K$ mentioned above, it is again clear that the
composition factors of $\Gamma$ and $G$ coincide.  It is therefore
natural to ask whether, in general, the existence of a Hopf-Galois
structure of type $G$ on a Galois extension with group $\Gamma$
necessarily forces $\Gamma$ and $G$ to have the same composition
factors.  We will answer this question in the negative by proving the
following result:

\begin{theorem} \label{index}
Let $G$ be a finite nonabelian simple group containing a subgroup $H$
of prime-power index, $|G:H|=p^a$ for $p$ prime and $a \geq 1$. Then
there exists a subgroup $J$ of $G$ of order $p^a$ such that any Galois
extension of fields with Galois group $\Gamma = H \times J$ admits a
Hopf-Galois structure of type $G$.
\end{theorem}

The proof of Theorem \ref{index} uses Guralnick's determination
\cite{Gu} of the nonabelian simple groups with a subgroup of
prime-power index (which depends on the classification of finite simple groups, and is recalled as Theorem \ref{guralnick} in \S\ref{differ} below),
together with the method of constructing Hopf-Galois structures via
fixed-point free pairs of homomorphism, introduced in \cite{CCo}
in the case $G=\Gamma$ and extended to the case $G \neq \Gamma$ in
\cite{BC}.

In particular, there exist finite soluble groups $\Gamma$ such that
every Galois extension with group $\Gamma$ admits a Hopf-Galois
structure of insoluble type. We collect together some examples of this
phenomenon in Corollary \ref{sol-insol} below, where we use the
following notation. As usual, $C_n$, $S_n$, $A_n$ and $D_{2n}$ denote
respectively the cyclic group of order $n$, the symmetric and
alternating groups of degree $n$, and the dihedral group of order
$2n$. For a prime $p$, we write $\PSL_m(p)$ for the projective
special linear group of dimension $m$ over the field $\F_p$ of $p$
elements. When $p$ is odd, we also write $F_{\frac{1}{2}p(p-1)}$ for the unique
Frobenius group of order $\frac{1}{2}p(p-1)$ that has a faithful
permutation representation of degree $p$, namely the semidirect
product $C_p \rtimes C_{\frac{1}{2}(p-1)}$ in which the second factor
acts faithfully on the first.

\begin{corollary} \label{sol-insol}
\begin{itemize}
\item[(i)] Any Galois extension of degree 60 with (soluble) Galois
  group $\Gamma = A_4 \times C_5$ admits a Hopf-Galois structure of
  (simple) type $A_5$.  
\item[(ii)] Any Galois extension of degree 168 with (soluble) Galois
  group $\Gamma = S_4 \times C_7$ admits a Hopf-Galois structure of 
  (simple) type $\PSL_3(2)\cong \PSL_2(7)$. 
\item[(iii)] if $p=2^e-1 \geq 7$ is a Mersenne prime, then any Galois
  extension of degree $\frac{1}{2}(p-1)p(p+1)$ with (soluble) Galois
  group $\Gamma = F_{\frac{1}{2}p(p-1)} \times D_{2^e}$ admits
  a Hopf-Galois structure of (simple) type $\PSL_2(p)$. 
\end{itemize}
\end{corollary}

We do not have any examples where an extension with insoluble Galois
group $\Gamma$ admits a Hopf-Galois structure of soluble type $G$. 

\section{Preliminaries on Hopf-Galois structures}

Let $L/K$ be a field extension of finite degree $n$, and let $H$ be a $K$-Hopf
algebra with comultiplication $\Delta \colon H \lra H \otimes_K H$,
$\Delta(h)=\sum_{(h)} h_{(1)} \otimes h_{(2)}$, and augmentation
(counit) $\epsilon \colon H \lra K$. We say that $L$ is an $H$-module
algebra if $H$ acts on $L$ so that $h \cdot (xy) = \sum_{(h)}
(h_{(1)}\cdot x) (h_{(2)}\cdot y)$ and $h \cdot k = \epsilon(h)k$ for
all $h \in H$, $x$, $y \in L$ and $k \in K$. We say that $L$ is an
$H$-Galois extension of $K$ if, furthermore, the $K$-linear map
$\theta \colon L \otimes_K L \lra \Hom_K(H,L)$, given by $\theta(x
\otimes y)(h) = x (h \cdot y)$ for $x$, $y \in L$ and $h \in H$,
is bijective. We then also say that $H$ endows $L$ with a Hopf-Galois
structure.

When $L/K$ is separable, Greither and Pareigis \cite{GP} gave the
following characterization of the Hopf-Galois structures on $L/K$. Let
$E$ be the normal closure of $L/K$ and consider the Galois groups
$\Gamma=\Gal(E/K)$ and $\Gamma'=\Gal(E/L)$. Let $\Gamma/\Gamma'$
denote the set of left cosets $\gamma \Gamma'$ of $\Gamma'$ in
$\Gamma$, and let $\Perm(\Gamma/\Gamma')$ be the group of permutations
of this set. Thus $\Perm(\Gamma/\Gamma')$ is isomorphic to the
symmetric group $S_n$. The left translations by elements of $\Gamma$
form a subgroup $\lambda(\Gamma)$ of $\Perm(\Gamma/\Gamma')$, which we
identify with $\Gamma$. The Hopf-Galois structures on $L/K$, up to
isomorphism, correspond bijectively with regular subgroups $G$ of
$\Perm(\Gamma/\Gamma')$ which are normalized by
$\lambda(\Gamma)$. (A permutation group acting on a set
$X$ is said to be regular if it is transitive and the stabilizer of any point in
$X$ is the identity.) In the Hopf-Galois structure corresponding to
such a subgroup $G$, the Hopf algebra acting on $L$ is
$E[G]^{\Gamma}$, the fixed point algebra of the group algebra $E[G]$
under the action of $\Gamma$ simultaneously on $E$ by field
automorphisms and on $G$ by conjugation by left translations inside
$\Perm(\Gamma/\Gamma')$. We refer to the isomorphism type of $G$ as
the {\em type} of the Hopf-Galois structure.

If $L/K$ is also normal, which will always be the case in this paper,
then $E=L$, $\Gamma'$ is trivial, and $\Gamma = \Gal(L/K)$. The result
of Greither and Pareigis then simplifies to the statement given in the
introduction. In particular, $\Gamma$ and $G$ are finite groups of the
same order $n$.  Moreover, if $\Gamma$ normalizes $G$ in $\Perm(\Gamma)$
then we can view $\Gamma$ as being contained in the holomorph $\Hol(G)
= G \rtimes \Aut(G)$ of $G$, which is usually a much smaller group
than $\Perm(\Gamma)$. Thus, given an abstract group $G$ of the same
order as $\Gamma$, there exists a Hopf-Galois structure on $L/K$ of
type $G$ if and only if there is an embedding $\beta \colon \Gamma
\lra \Hol(G)$ with regular image. We shall call such a $\beta$ a
regular embedding. Then the number of Hopf-Galois structures of type
$G$ is the number of equivalence classes of regular embeddings, where two
embeddings are equivalent if and only if they are conjugate in
$\Hol(G)$ by an element of $\Aut(G)$; see \cite{unique} or
\cite[\S7]{Ch00}.

It is convenient to write elements of $\Hol(G)$ in the form
$[g,\alpha]$ with $g \in G$ and $\alpha \in \Aut(G)$. The
multiplication in $\Hol(G)$ is then given by
\begin{equation} \label{hol-mult}
   [g_1, \alpha_1] [g_2 ,\alpha_2] = [g_1 \alpha_1(g_2), \alpha_1 \alpha_2 ]. 
\end{equation}

\section{Nilpotent Hopf-Galois structures}

In this section, we prove Theorem \ref{nilp-sol}. Thus we suppose that
$L/K$ is a finite Galois extension of fields with group
$\Gal(L/K)=\Gamma$, and that $L/K$ admits a Hopf-Galois structure of
type $G$ for some nilpotent group $G$. We then have a regular
embedding $\beta : \Gamma \lra \Hol(G)$. Moreover, being nilpotent,
$G$ can be written as the direct product
 $$ G=\prod_p G_p  $$ 
over all primes $p$, where $G_p$ denotes the
(unique) Sylow $p$-subgroup of $G$ \cite[5.2.4]{Rob}. Our task is to
show that the group $\Gamma$ is soluble. 
 
If $J$ is a finite group of order $p^r m$, where $p$ is prime and $p
\nmid m$, then a subgroup of $J$ of order $m$ is called a Hall
$p'$-subgroup of $J$. We will use P.~Hall's theorem \cite[9.1.8]{Rob} that
if $J$ has a Hall $p'$-subgroup for every prime $p$, then $J$ is
soluble.

For each $p$, let 
$$ H_p = \prod_{q \neq p} G_q. $$ 
Then $H_p$ is a Hall $p'$-subgroup of $G$. Moreover, $H_p$ is a
characteristic subgroup of $G$, as it consists precisely of the
elements of order prime to $p$.

Define $\Delta_p=\{ \gamma \in \Gamma \: : \:  \beta(\gamma)\cdot e_G \in
H_p \}$, where $e_G$ is the identity element of $G$.  Since
$\beta(\Gamma)$ is regular on $G$, it is clear that $\Delta_p$ is a {\em
  subset} of $\Gamma$ of size $|H_p|$. However, since $H_p$ is
characteristic in $G$, more is true:

\begin{lemma} 
$\Delta_p$ is a subgroup of $\Gamma$.
\end{lemma}
\begin{pf}
As $\Delta_p$ is nonempty and $\Gamma$ is finite, it suffices to check
that $\Delta_p$ is closed under multiplication. So let $\sigma_1$,
$\sigma_2 \in \Delta_p$, and let $h_i=\beta(\sigma_i) \cdot e_G$ for
$i=1$, $2$. Then $h_i \in H_p$, and, in the notation of
(\ref{hol-mult}), we have $\beta(\sigma_i)=[h_i,\alpha_i]$ for some
$\alpha_i\in \Aut(G)$. Then $\beta(\sigma_1 \sigma_2) = [h_1 \alpha_1(h_2),
\alpha_1 \alpha_2]$, so that $\beta(\sigma_1 \sigma_2) \cdot e_G = h_1
\alpha(h_2)$. Since $H_p$ is a characteristic subgroup of $G$, we have
$\alpha_1(h_2) \in H_p$ and hence $\beta(\sigma_1 \sigma_2) \cdot e_G\in
H_p$. Thus $\sigma_1 \sigma_2 \in \Delta_p$, as required.
\end{pf}

As $|\Gamma|=|G|$ and $|\Delta_p|=|H_p|$, it follows that $\Delta_p$ is
a Hall $p'$-subgroup of $\Gamma$. Thus $\Gamma$ contains a Hall
$p'$-subgroup $\Delta_p$ for each prime $p$. Hence, by Hall's theorem, $\Gamma$ is soluble. This completes the proof of Theorem \ref{nilp-sol}.

\section{Abelian extensions}

In this section and the next, we prove Theorem \ref{ab-sol}. We first
reduce the proof of Theorem \ref{ab-sol} to that of a statement,
Theorem \ref{a-ineq}, about simple groups. 
Theorem \ref{a-ineq} will then be proved in \S\ref{simple-gps}, using
a theorem of Vdovin \cite{V} (Theorem \ref{V-thm} in \S\ref{simple-gps} below) which is a consequence
of the classification of the finite simple groups. 

Before stating Theorem \ref{a-ineq}, we introduce some notation.

\begin{defn}
Let $G$ be a finite group. Then 
$$ a(G) = \max\{ \; |A| \; : \; A \mbox{ is an abelian subgroup of } G
\}. $$
\end{defn}

We note another result from Vdovin's paper, and record some obvious
properties of $a(G)$. 

\begin{proposition} \label{aSn}
For the symmetric group $S_m$, we have $a(S_m) \leq 3^{m/3}$.
\end{proposition}
\begin{pf}
From \cite[Theorem 1.1]{V}, we have $a(S_{3k})=3^k$, $a(S_{3k+1})=4
\cdot 3^{k-1}$ and $a(S_{3k+2})=2 \cdot 3^k$. Hence the stated
inequality holds in all cases.
\end{pf}

\begin{proposition} \label{a-prop}
Let $G$ be a finite group.
\begin{itemize}
\item[(i)] If $H$ is a subgroup of $G$ then $a(H) \leq a(G)$.
\item[(ii)] If $N$ is a normal subgroup of $G$ then $a(G) \leq
  a(N)a(G/N)$.
\item[(iii)] If $G=H \times J$ then $a(G)=a(H) a(J)$. In particular $a(H^m)=a(H)^m$.
\end{itemize}
\end{proposition}
\begin{pf}

\noindent (i) Any abelian subgroup of $H$ is certainly an abelian subgroup of $G$.

\noindent (ii) Let $A$ be an abelian subgroup of $G$, and let $A_1$
  (respectively, $A_2$) be the kernel (respectively, image) of the
  composite homomorphism $A \hookrightarrow G \twoheadrightarrow
  G/N$. Then $|A| = |A_1| \, |A_2| \leq a(N) a(G/N)$. 

\noindent (iii) Let $G=H \times J$. By (ii), $a(G) \leq
a(H) a(J)$. But if $A \subseteq H$, $B \subseteq J$ are abelian
subgroups with $|A|=a(H)$, $|B|=a(J)$, then $A \times B$ is an abelian
subgroup of $G$ of order $a(H) a(J)$. This gives the first assertion, and the second follows by induction.
\end{pf}
 
We will deduce Theorem \ref{ab-sol} from the following statement:

\begin{theorem} \label{a-ineq}
Let $T$ be a finite nonabelian simple group. Then   
\begin{equation} \label{a-ineq-1}
      3^{1/3} a(T) a(\Aut(T)) < |T|. 
\end{equation}
\end{theorem}

\begin{pf-of}{Theorem \ref{ab-sol} (assuming Theorem \ref{a-ineq})}

Let $\Gamma$ be a finite abelian group. We need to show that if there
is a regular embedding $\beta \colon \Gamma \lra \Hol(G)$ for some
group $G$, then $G$ must be soluble. 

We first treat the special case where $G$ is characteristically simple. Thus, by
\cite[3.3.15]{Rob}, $G$ is the direct product $T^m$  
for some simple group $T$ and some $m \geq 1$. We claim that
the existence of a regular embedding $\Gamma \lra \Hol(T^m)$, with
$\Gamma$ abelian, forces $T$ to be abelian. Then $G$ is abelian, and
hence soluble, as required. 

Suppose that $T$ is a nonabelian simple group. Then, by \cite[Lemma
  3.2]{simple}, $\Aut(T^m)$ is the wreath product
$$ \Aut(T^m) = \Aut(T) \mathrm{\wr\ } S_m = (\Aut(T)^m) \rtimes S_m , $$
where the symmetric group $S_m$ permutes the $m$ factors.
We therefore have a regular embedding of the abelian group $\Gamma$ in 
$$ \Hol(T^m) = T^m \rtimes ( \Aut(T)^m \rtimes S_m). $$
We break $\Hol(T^m)$ into the sequence of quotients
$$ H_1= \frac{\Hol(T^m)}{\Hol(T)^m} \cong \frac{T^m \rtimes(\Aut(T)^m
  \rtimes S_m)}{(T \rtimes \Aut(T))^m} \cong S_m, $$
$$ H_2 = \frac{\Hol(T)^m}{T^m} \cong \Aut(T)^m; $$
$$ H_3 = T^m. $$

As $|T|^m=|G|=|\Gamma|=|\beta(\Gamma)|$, we may apply Proposition
\ref{a-prop}(ii) twice to get $|T|^m  \leq a(H_1) a(H_2) a(H_3)$.
Now $a(H_1) \leq 3^{m/3}$ by Proposition \ref{aSn}, and 
$a(H_2)=a(\Aut(T)^m)=a(\Aut(T))^m$ and $a(H_3)=a(T)^m$ by Proposition
\ref{a-prop}(iii). Thus we have
$$ |T|^m \leq 3^{m/3} a(\Aut(T))^m a(T)^m, $$
contradicting Theorem \ref{a-ineq}. Hence $T$ is
abelian, as claimed, and Theorem \ref{ab-sol} holds when $G$ is
characteristically simple.

We now prove the general case by induction on $|G|=|\Gamma|$, taking
the case just considered as the base of the induction. 
If $G$ is not characteristically simple then it has a  
nontrivial proper characteristic subgroup $H$, and by 
\cite[Proposition 3.1]{simple},  $\beta$ induces a
homomorphism
$$ \overline{\beta} : \Gamma \lra \Hol(G/H) $$
whose image is transitive on $G/H$. This image is also abelian, and is
therefore regular on $G/H$. Let
$\Sigma=\ker(\overline{\beta})$. Then $|\Sigma|=|H|$ and the abelian
group $\Sigma$ acts regularly on $H$, so $\beta$ restricts to a regular embedding $\Sigma \lra \Hol(H)$. As $|H|$, $|G/H| < |G|$, it
follows from the induction hypothesis that $H$ and $G/H$ are soluble,
and hence so is $G$.
\end{pf-of}

\section{Proof of Theorem \ref{a-ineq}} \label{simple-gps} 

To complete the proof of Theorem \ref{ab-sol}, we must prove Theorem
\ref{a-ineq}. We do so using the classification of finite simple
groups. Our main reference for the necessary facts about the simple
groups is the book \cite{W}. We also use \cite{Gor}. In brief, the
classification states that every finite
nonabelian simple group is either an alternating group $A_n$, $n \geq
5$, a (classical or exceptional) group of Lie type, or one of the $26$
sporadic simple groups. We refer the reader to \cite[p.~3]{W}
for a more detailed statement.

Using the classification, Vdovin \cite[Theorem A]{V} proved the
following result:

\begin{theorem}[Vdovin] \label{V-thm}
Let $T$ be a finite nonabelian simple group which is not of the form
$\PSL_2(q)$. Then $a(T)^3 < |T|$.
\end{theorem}

\begin{remark}
Some of the groups $\PSL_2(q)$ appear in the classification in another
guise, so that some families in the classification other than $\PSL_2(q)$
contain groups for which the conclusion of Vdovin's theorem does not
hold. Thus, for the alternating group $A_5$ of order $60$, we have $a(A_5)=5 >
60^{1/3}$, but this does not contradict Vdovin's theorem since $A_5
\cong \PSL_2(4)\cong \PSL_2(5)$. 
\end{remark}

For any nonabelian simple group $T$, the group $\Aut(T)$ contains the
subgroup $\Inn(T) \cong T$ of inner automorphisms, and we write 
$$ \Out(T) = \frac{\Aut(T)}{\Inn(T)} $$ 
for the group of outer automorphisms. Using Proposition
\ref{a-prop}(ii), we then have
\begin{equation} \label{Aut}
   a(\Aut(T)) \leq a(\Inn(T)) a(\Out(T)) \leq a(T) |\Out(T)|. 
\end{equation}
A famous consequence of the
classification is the proof of the Schreier Conjecture, which asserts that
$\Out(T)$ is always soluble. Of greater relevance for us is the fact
that $\Out(T)$ is very small relative to $T$. We will use this fact in conjunction with the following result.

\begin{proposition} \label{alt-ineq}
Let $T$ be a nonabelian simple group which is not of the form
$\PSL_2(q)$. If the inequality 
\begin{equation} \label{a-ineq-2}
  3 |\Out(T)|^3 < |T| 
\end{equation}
holds, then the conclusion (\ref{a-ineq-1}) of Theorem \ref{a-ineq} holds for $T$.
\end{proposition}
\begin{pf}
Since $a(T)<|T|^{1/3}$ by Theorem \ref{V-thm}, it follows from
(\ref{Aut}) and (\ref{a-ineq-2}) that 
\begin{eqnarray*}
 3^{1/3} a(T) a(\Aut(T)) & \leq & 3^{1/3} a(T)^2 |\Out(T)| \\
                        & < & a(T)^2 |T|^{1/3} \\
                        & < & |T| . 
\end{eqnarray*}
Thus (\ref{a-ineq-1}) holds.
\end{pf}

By the classification of finite simple groups, together with
Proposition \ref{alt-ineq}, the proof of Theorem \ref{a-ineq} is
reduced to the following five lemmas, whose proofs will occupy the
rest of this section.

\begin{lemma} \label{alt}
The conclusion (\ref{a-ineq-1}) of Theorem \ref{a-ineq} holds for each
alternating group $T=A_n$, $n \geq 5$.
\end{lemma}

\begin{lemma} \label{sporadic}
The inequality (\ref{a-ineq-2}) of Proposition \ref{alt-ineq} holds
for each sporadic simple group $T$.
\end{lemma}

\begin{lemma}  \label{PSL2}
The conclusion (\ref{a-ineq-1}) of Theorem \ref{a-ineq} holds for each
simple group $T$ of the form $\PSL_2(q)$.
\end{lemma}

\begin{lemma} \label{classical}
The inequality (\ref{a-ineq-2}) of Proposition \ref{alt-ineq} holds
for each simple group $T$ which is a classical group of Lie type but
is not of the form $\PSL_2(q)$.
\end{lemma}

\begin{lemma} \label{exceptional}
The inequality (\ref{a-ineq-2}) of Proposition \ref{alt-ineq} holds
for each simple group $T$ which is an exceptional group of Lie type.
\end{lemma}

The first two cases are easily handled. \\

\begin{pf-of}{Lemma \ref{alt}}
We treat all the simple alternating groups $T=A_n$, $n\geq 5$ here,
even though some of them are of the form $\PSL_2(q)$. In all cases, we
have $a(A_n) \leq a(S_n) \leq 3^{n/3}$ by Propositions \ref{a-prop}(i)
and \ref{aSn}. For $T=A_n$ with $n \neq 6$, we have $|\Out(A_n)|=2$
and $\Aut(A_n) \cong S_n$, whereas $|\Out(A_6)|=4$ and $\Aut(A_6)$ has
a normal subgroup of index $2$ isomorphic to $S_6$ \cite[pp.~18,
  19]{W}. Thus for $n \neq 6$, we have $3^{1/3} a(T) a(\Aut(T) \leq
3^{(2n+1)/3}$ and (\ref{a-ineq-1}) will hold provided that
$$ 3^{(2n+1)/3} < \textstyle{\frac{1}{2}} n!. $$
But this inequality holds for $n=5$ and hence, by induction, for all
$n \geq 5$. In the exceptional case $T=A_6$, we have $a(\Aut(T))) \leq
2 a(S_6) \leq 18$, so $3^{1/3} a(T)a(\Aut(T)) \leq 2 \cdot 3^{13/3} <
360=|T|$, and again (\ref{a-ineq-1}) holds. 
\end{pf-of}

\begin{pf-of}{Lemma \ref{sporadic}}
Let $T$ be a sporadic simple group. From
\cite[p.~304]{Gor}, we have
$|\Out(T)| \leq 2$.
Thus $3|\Out(T)|^3 \leq 24 < |T|$, so (\ref{a-ineq-2}) holds. 
\end{pf-of}

Before giving the proofs of Lemmas \ref{PSL2}--\ref{exceptional}, 
we recall some
general properties of the finite simple groups of Lie type and their outer automorphisms, and we summarize the facts we will need about the various families of groups in Tables \ref{lie-class}--\ref{lie-except-autom} below. The information in these tables is
taken from \cite[Chapters 3, 4]{W} and \cite[p.~135]{Gor}.  The outer
automorphisms for most of the groups of Lie type are also described in
\cite{St}.

The groups in each family are indexed by a prime-power parameter $q$, and
we will always write $q=p^e$ with $p$ prime. Each simple group $T$ of
Lie type is obtained as the central quotient of some group of
matrices over a finite field. Following \cite[p.~135]{Gor}, and
changing notation from the previous sections, we denote this group by
$G$, and we write $d$ for the order of its center. We therefore have $|T| = |G|/d$.  For example, if $T$ is
the projective special linear group $\PSL_n(q)$ then the corresponding
group $G$ is $\SL_n(q)$ whose center has order $d = (n,q-1)$. We list the
classical simple groups $T$ of Lie type in Table \ref{lie-class}, together
with the corresponding groups $G$ and their orders. The value of $d$ is
shown in Table \ref{lie-class-autom}, along with the quantities
$\epsilon$ and $g$ which will be explained below. The
exceptional simple Lie groups $T$, and the orders of the corresponding
groups $G$, are given in Table \ref{lie-except}, with the values of
$d$, $\epsilon$ and $g$ for each group shown in Table
\ref{lie-except-autom}. We use the notation of \cite{W} for the groups
$T$, but the notation of \cite{Gor} for the groups $G$. The notation
for the groups $G$ is derived from the standard labelling of the
associated Dynkin diagrams.  For example, the group $\SL_n(q)$ is
denoted $A_{n-1}(q)$ as it corresponds to the Dynkin diagram
$A_{n-1}$. 

We have the generic isomorphism $\PO_5(q) \cong \PSp_4(q)$
\cite[p.~96]{W}. Following \cite{W}, we omit the groups $\PO_5(q)$
from the classification (whereas \cite{Gor} omits the groups
$\PSp_4(q)$ instead). Also, $\PO_{2n+1}(q) \cong \PSp_{2n}(q)$ when
$q=2^e$, so the groups $\PO_{2n+1}(q)$ for $q$ even could be omitted.

For each simple group $T$ of Lie type, any outer automorphism of $T$
may be written as a product of a diagonal automorphism, a field
automorphism and a graph automorphism.  The diagonal automorphisms
arise from conjugation of $G$ by elements of a larger matrix group in which $G$ is normal. For example, the diagonal automorphisms of $\PSL_n(q)$ are induced by the automorphisms of $\SL_n(q)$ arising from conjugation by elements of $\GL_n(q)$. In all cases, the number of such diagonal automorphisms is the quantity $d$ described above. The field automorphisms are induced by automorphisms of the
underlying finite field, and therefore form a cyclic group. We write
$\epsilon$ for the number of field automorphisms. In most cases,
$\epsilon = e$ as we are working with matrices over the field $\F_q$
of $q=p^e$ elements.  The exceptions are that $\epsilon=2e$ for the
families $\PSU_n(q)$, $\PO^-_{2n}(q)$ and ${}^2E_6(q)$, since these groups
are obtained from groups of matrices over $\F_{q^2}$, and
$\epsilon=3e$ for ${}^3D_6(q)$ where the matrices are over
$\F_{q^3}$. Finally, the graph automorphisms arise from automorphisms
of the Dynkin diagram where, for groups in characteristic $p$, the
automorphism does not need to preserve the direction of the arrow on
an edge of multiplicity $p$. We write $g$ for the number of graph
automorphisms. We then have $|\Out(T)|=d \epsilon g$
\cite[pp.~303, 304]{Gor}. Thus Tables \ref{lie-class-autom} and
\ref{lie-except-autom} enable us to find $|\Out(T)|$ in all cases. 

\begin{table} 
\centerline{
\begin{tabular}{|c|c|c|c|}
\hline
  $T$ & restrictions & $G$ & $|G|$   \\
\hline
$\PSL_n(q)$ & $n \geq 2$ & $A_{n-1}(q)$ & 
 $q^{\frac{1}{2}n(n-1)} \prod_{i=2}^n (q^i-1)$ \\
            & $(n,q) \neq (2,2)$ & &   \\ 
       & $(n,q) \neq (2,3)$ & &   \\ 
\hline
$\PSU_n(q)$ & $n \geq 3$ & ${}^2A_{n-1}(q)$ & 
  $ q^{\frac{1}{2}n(n-1)}\prod_{i=2}^n \bigl(q^i-(-1)^i\bigr)$ \\
  & $(n,q) \neq (3,2)$ & &  \\
\hline
$\PSp_{2n}(q)$ & $n \geq 2$ & $C_n(q)$ & 
 $q^{n^2} \prod_{i=1}^n (q^{2i}-1)$  \\
 & $(n,q) \neq (2,2)$ & &  \\ 
\hline
$\PO_{2n+1}(q)$ & $n \geq 3$ & $B_n(q)$ & 
 $q^{n^2} \prod_{i=1}^n (q^{2i}-1)$    \\  
\hline
$\PO^+_{2n}(q)$ & $n \geq 4$ & $D_n(q)$ & 
 $q^{n(n-1)} (q^n-1)\prod_{i=1}^{n-1} (q^{2i}-1)$   \\
\hline
$\PO^-_{2n}(q)$ & $n \geq 4$ & ${}^2D_n(q)$ & 
 $q^{n(n-1)} (q^n+1)\prod_{i=1}^{n-1} (q^{2i}-1)$   \\  
\hline
\end{tabular} 
} 
\caption{Classical simple groups of Lie type} \label{lie-class}  
\end{table}

\begin{table} 
\centerline{
\begin{tabular}{|c|c|c|c|c|}
\hline
  $T$ & restrictions & $d$ & $\epsilon$ & $g$ \\
\hline
$\PSL_n(q)$ & $n \geq 2$  & $(n,q-1)$ & $e$ &
$1$ if $n=2$ \hfill \\
            & $(n,q) \neq (2,2)$ & & &  $2$ if $n>2$ \hfill \\ 
       & $(n,q) \neq (2,3)$ & & &   \\ 
\hline
$\PSU_n(q)$ & $n \geq 3$ &
$(n,q+1)$ & $2e$ & $1$ \\
  & $(n,q) \neq (3,2)$ & & & \\
\hline
$\PSp_{2n}(q)$ & $n \geq 2$ & $(2,q-1)$ & $e$ & $2$ if $n=2$ and $q$
is even  \\
 & $(n,q) \neq (4,2)$ & & & $1$ otherwise \\ 
\hline
$\PO_{2n+1}(q)$ & $n \geq 3$ & $(2,q-1)$ & $e$ & $1$  \\  
\hline
$\PO^+_{2n}(q)$ & $n \geq 4$ &  $(4,q^n-1)$ & $e$
  &  $6$ if $n=4$ \\
 & & &  & $2$ otherwise \\
\hline
$\PO^-_{2n}(q)$ & $n \geq 4$ & $(4,q^n+1)$ & $2e$ & $1$  \\  
\hline
\end{tabular}
}   
\caption{Automorphisms of classical simple groups of Lie type}
 \label{lie-class-autom}  
\end{table}

\begin{table}
\centerline{
\begin{tabular}{|c|c|}
\hline
  $T$ & $|G|$   \\
\hline
$G_2(q)$  & $q^6(q^6-1)(q^2-1)$   \\
  ($q \geq 3$)  &  \\
\hline 
$F_4(q)$ &   $q^{24}(q^{12}-1)(q^8-1)(q^6-1)(q^2-1)$   \\
\hline
$E_6(q)$ &   $q^{36}(q^{12}-1)(q^9-1)(q^8-1)(q^6-1)(q^5-1)(q^2-1)$   \\
\hline
${}^2E_6(q)$  & $q^{36}(q^{12}-1)(q^9+1)(q^8-1)(q^6-1)(q^5+1)(q^2-1)$  \\
\hline
${}^3D_4(q)$   & $q^{12}(q^8+q^4+1)(q^6-1)(q^2-1)$ \\
\hline
$E_7(q)$ &  $q^{63}(q^{18}-1)(q^{14}-1)(q^{12}-1)(q^{10}-1)(q^8-1)(q^6-1)(q^2-1)$  \\
\hline
$E_8(q)$ &    $q^{120} (q^{30}-1)(q^{24}-1)(q^{20}-1)(q^{18}-1) (q^{14}-1)(q^{12}-1)(q^8-1) (q^2-1) $  \\
\hline
${}^2B_2(q)$ & $q^2(q^2+1)(q-1)$  \\
 ($q=2^{2n+1}$ &  \\
$n \geq 1$) &   \\
\hline
${}^2G_2(q)$ & $ q^3 (q^3+1) (q-1)$    \\
($q=3^{2n+1}$ &   \\
$n \geq 1$) &  \\
\hline
${}^2F_4(q)$ & $q^{12}(q^6+1)(q^4-1)(q^3+1)(q-1)$  \\ 
($q=2^{2n+1}$   &  \\
$n \geq 1$) &    \\
\hline
${}^2F_4(2)'$ &  $q^{12}(q^6+1)(q^4-1)(q^3+1)(q-1)$ with $q=2$  \\ 
\hline
\end{tabular} 
}  
\caption{Exceptional simple groups of Lie type}  \label{lie-except}    
\end{table}

\begin{table} 
\centerline{ 
\begin{tabular}{|c|c|c|c|}
\hline
  $T$ & $d$ & $\epsilon$ & $g$ \\
\hline
$G_2(q)$  & $1$ & $e$ & $2$ if $q=3^e$ \, \hfill \\
          &     &     & $1$ otherwise \hfill \\
\hline 
$F_4(q)$ & $1$ & $e$ & $2$ if $q=2^e$ \, \hfill \\
         &     &     & $1$ otherwise \hfill \\
\hline
$E_6(q)$ & $(3,q-1)$ & $e$ & $2$ \\
\hline
${}^2E_6(q)$  & $(3,q+1)$ & $2e$ & $1$  \\
\hline
${}^3D_4(q)$ & $1$ & $3e$ & $1$  \\
\hline
$E_7(q)$ & $(2,q-1)$ & $e$ &  $1$ \\
\hline
$E_8(q)$ & $1$ & $e$ & $1$  \\
\hline
${}^2B_2(q)$ & $1$ & $e$ & $1$ \\
\hline
${}^2G_2(q)$ & $1$ & $e$ & $1$ \\
\hline
${}^2F_4(q)$  & $1$ & $e$ & $1$  \\ 
\hline
${}^2F_4(2)'$ & $2$   & $1$ & $1$ \\
\hline
\end{tabular}
}  
\caption{Automorphisms of exceptional simple groups of Lie type}
  \label{lie-except-autom}     
\end{table}

We now give the proofs of Lemmas \ref{PSL2}--\ref{exceptional}.
\\

\begin{pf-of}{Lemma \ref{PSL2}} 
Let $T=\PSL_2(q)$. In this case, Theorem \ref{V-thm} and Proposition
\ref{alt-ineq} do not apply. Note that, for $T$ to be simple, we require
$q \geq 4$.  Moreover, as $\PSL_2(4) \cong \PSL_2(5) \cong A_5$ and
$\PSL_2(9) \cong A_6$ (see \cite[p.~3]{W}), the cases $q=4$, $5$ and
$9$ follow from Lemma \ref{alt}. We verify (\ref{a-ineq-1}) for the
remaining values of $q$.

We have 
$$ a(T) = \begin{cases} q+1 & \mbox{if } q=2^e, \\
                        q   & \mbox{if } q \mbox{ is odd};
           \end{cases}   $$
see \cite[Theorem 3.1]{V} or the full list of subgroups of $T$ in
\cite[8.27]{Hup}. Also, $|\Out(T)|=de$ with $d=(q-1,2)$ by Table
\ref{lie-class-autom} (or \cite[Theorem 3.2, p.~50]{W}). 

If $q=2^e$ with $e \geq 3$, then $d=1$, $|T|=q(q^2-1)$, $|\Out(T)|= e$ and
$a(T)=q+1$. Thus $a(\Aut(T)) \leq (q+1)e$. So it suffices to show that 
$3^{1/3} (q+1)^2 e < q(q^2-1)$, which will follow from $3^{1/3} e<
q-2$. The last inequality holds for $e=3$ (so $q=8$), and hence for
all $e \geq 3$. 

For odd $q>3$, we have $d=2$, $|\Out(T)|=2e$, $|T|=\frac{1}{2}q(q^2-1)$ and
$a(T)=q$. It therefore suffices to show that
$$ 3^{1/3} \cdot 2 e q^2 < \textstyle{\frac{1}{2}} q (q^2-1), $$
which will follow from 
$$ 3^{1/3} \cdot 4 e < q-1. $$
This holds for $q=p^e$ if $ p \geq 7$, $e \geq 1$, if $p=5$, $e \geq
2$ and $p=3$, $e \geq 3$. Hence (\ref{a-ineq}) holds in all cases.
\end{pf-of}
 
In the proofs of the remaining two lemmas, we will frequently use the fact
that, for any prime power $q=p^e$, $e \geq 1$, we have $e^3 \leq
\frac{1}{2}q^2$. We also note that, since $|T|=|G|/d$, the inequality 
(\ref{a-ineq-2}) to be proved may be rewritten as
\begin{equation} \label{a-ineq-3}
  3 d |\Out(T)|^3 < |G|.
\end{equation}

\begin{pf-of}{Lemma \ref{classical}}
Let $T$ be a classical group of Lie type
which is not of the form $\PSL_2(q)$. 
We consider each family in Table \ref{lie-class}. 

For $T=\PSL_n(q)$ with $n \geq 3$, we have $d \leq q-1$, $g=2$ so that
$|\Out(T)| \leq 2(q-1)e$. Thus
$$ 3d |\Out(T)|^3 \leq 24 (q-1)^4 e^3 \leq 12 (q-1)^4 q^2. $$
If $n \geq 4$ then $|G| > q^{12}$ so $3d |\Out(T)|^3 < 2^4 q^6 < |G|$.
If $n=3$ then $|G|>q^7(q-1)$, so that $3d |\Out(T)|^3  <
12q^5(q-1) < |G|$ provided that $q \geq 4$. For $n=3$ and $q=2$ or $3$, we
have $d=1$ so that $3d |\Out(T)|^3 = 24 e^3 \leq 12 q^2 < q^6 < |G|$.
Thus (\ref{a-ineq-3}) holds for all the simple groups $\PSL_n(q)$. 

For $T=\PSU_n(q)$, we have $d \leq q+1 \leq \frac{3}{2}q$ and
$$  3d|\Out(T)|^3 \leq 3 (q+1)^4 (2e)^3.  $$ 
So if $n \geq 4$ then 
$$ 3d|\Out(T)|^3 \leq 3 \left( \textstyle{\frac{3}{2}}q \right)^4 
 (2e)^3 < 2^6 q^6 \leq q^{12} < |G|. $$ 
Now let $n=3$ (so $q \geq 3$). Then $|G| > q^7(q-1)$, and, since $d \leq q+1 \leq \frac{4}{3}q$, we have
$$ 3d|\Out(T)|^3 \leq 3 \left( \textstyle{\frac{4}{3}} q \right)^4
(2e)^3 < 40 q^6. $$
Thus (\ref{a-ineq-3}) holds if $q \geq 7$, since then $q(q-1)>40$.
It remains to check the cases $q=3$, $4$ and $5$. We have $d=1$ for
$q=3$, $4$, and $d=3$ for $q=5$. Also $e=1$ for $q=3$, $5$ and $e=2$
for $q=4$. Hence $|\Out(T)| = 2de \leq 6$ for $q \leq 5$, so that 
$3d |\Out(T)|^3 \leq 9 \cdot 6^3 <3^7 <|G|$. 

For $T=\PSp_{2n}(q)$ or $\PO_{2n+1}(q)$, we have $d \leq 2$, 
$g \leq 2$ so
$$ 3d|\Out(T)|^3 \leq 3 \cdot 2^7 e^3 \leq 3 \cdot 2^6 q^2 \leq
\textstyle{\frac{3}{4}} q^{10} < |G|, $$
where the last inequality holds as the case $(n,q)=(2,2)$ does not
occur. 

For $T=\PO_{2n}^+(q)$ or $T=\PO_{2n}^-(q)$, we have $d \leq 4$, and in
both cases $|\Out(T)|=d \epsilon g \leq 24e$. Thus
$$  3d|\Out(T)|^3 \leq 12 \cdot (24e)^3 \leq 2^{10} \cdot 3^4 
q^2 < 2^{17}q^2 \leq q^{19} < |G|. $$
\end{pf-of}

\begin{pf-of}{Lemma \ref{exceptional}}
Let $T$ be an exceptional group $T$ of
Lie type. 
From Table \ref{lie-except-autom}, we have $d \leq 3$ and
$|\Out(T)| \leq 6e$ in all cases.  Thus it
suffices to prove that $9 (6e)^3 < |G|$. This inequality will hold if
\begin{equation} \label{except-bound}
   4 \cdot 3^5 q^2 < |G|
\end{equation} 
Now for $q \geq 2$, we have $4 \cdot 3^5 < q^{10}$, so
(\ref{except-bound}) will hold if $|G| > q^{12}$. From Table
\ref{lie-except}, this covers all cases except
${}^2B_2(q)$ and ${}^2G_2$. For these two cases, we have $d=1$ and
$|\Out(T)|=e$, so that we only need to show $3e^3<|G|$.  But as
$e^3<q^2$ we have $3e^3<q^3<|G|$ in both cases. Thus the 
inequality (\ref{a-ineq-3}) holds for all the exceptional groups of
Lie type. 
\end{pf-of}

This concludes the proof of Theorem \ref{a-ineq} and thus of Theorem \ref{ab-sol}.

\section{Cases where the composition factors differ} \label{differ}
 
In this section, we show that a Galois extension with group $\Gamma$ may
admit Hopf-Galois structures of type $G$, where the composition
factors of the groups $\Gamma$ and $G$ differ. In particular, we shall
prove Theorem \ref{index} and Corollary \ref{sol-insol}.
All our examples arise from the construction we shall give in Lemma
\ref{untangle}. 
 
Two subgroups $H$, $J$ in a finite group $G$ are said to be
complementary if $|H| \, |J| =|G|$ and $H \cap J =\{e_G\}$. (We do not
require either subgroup to be normal.) Then
each element of $G$ can be written uniquely in the form $hj$
with $h \in H$ and $j \in J$. If $G$ acts faithfully and transitively
as permutations of some set $X$, and $H$ is the stabilizer of a point
in $X$, then a subgroup $J$ of $G$ is complementary to $H$ if and only
if $J$ is regular on $X$. In particular, if $G$ contains a Hall
$p'$-subgroup $H$ for some prime $p$ then (taking $X$ to be the space
of left cosets of $H$), any Sylow $p$-subgroup of $G$ is complementary
to $H$. 

\begin{lemma} \label{untangle}
Suppose that $G$ contains a pair of complementary subgroups $H$, $J$.
Then any Galois extension of fields with Galois 
group $\Gamma = H \times J$ admits a Hopf-Galois structure of type
$G$.
\end{lemma}
\begin{pf}
It suffices to exhibit a regular embedding $\beta : H \times J \lra
\Hol(G)$. Using the notation of (\ref{hol-mult}), we claim that such
an embedding is given by the formula 
$\beta(h,j)=[hj^{-1},C(j)]$ for $h \in H$, $j \in J$ where, for any $g
\in G$, $C(g)\in \Aut(G)$ is conjugation by $g$, that is,
$C(g)(x)=gxg^{-1}$ for $x \in G$. More concisely,
$\beta(h,j)(x)=hj^{-1}(jxj^{-1})=hxj^{-1}$ for $x \in G$. 

We check that $\beta$ is indeed a homomorphism. If $h_1$, $h_2 \in
H$ and $j_1$, $j_2 \in J$, then, for each $x \in G$, we calculate
$$  \beta(h_1,j_1) \beta(h_2,j_2)(x)  = h_1(h_2xj_2^{-1})j_1^{-1}
          = (h_1h_2) x (j_1 j_2)^{-1} = \beta(h_1 h_2, j_1 j_2)(x), $$
so that $\beta(h_1,j_1) \beta(h_2,j_2)=\beta(h_1 h_2,j_1 j_2)$ as
required. The homomorphism $\beta$ is regular since each $x 
\in G$ can be written uniquely in the form $hj^{-1}$ with $h \in H$ and $j
\in J$.  
\end{pf}

\begin{remark}
Lemma \ref{untangle} is an application of the method of fixed-point
free pairs of homomorphisms; see \cite[\S2]{BC} (or, in the case that
$\Gamma=G$, \cite[\S4]{CCo}). For finite groups $\Gamma$, $G$ of the
same order, we say that homomorphisms
$$ \beta_1, \beta_2 \colon \Gamma \lra G $$
form a fixed-point free pair if $\beta_1(\sigma)=\beta_2(\sigma)$ only
for $\sigma=e_\Gamma$. We then have a regular embedding $\beta \colon
\Gamma \lra \Hol(G)$ given by $\beta(\sigma)=\lambda(\beta_1(\sigma))
\rho(\beta_2(\sigma)) = [\beta_1(\sigma)\beta_2(\sigma)^{-1},
  C(\beta_2(\sigma))]$, where $\lambda, \rho 
    \colon G \lra \Hol(G)$ are the regular left and right embeddings,
    $\lambda(\sigma)(\tau)= \sigma \tau$, $\rho(\sigma)(\tau)=\tau
    \sigma^{-1}$ for $\tau \in G$. In the proof of Lemma
    \ref{untangle}, we have $\beta_1, \beta_2 \colon H \times J \lra
    G$ with $\beta_1(h,j)=h$ and $\beta_2(h,j)=j$. 
 \end{remark}

Before proving Theorem \ref{index}, we consider the case of symmetric
and alternating groups.

\begin{example} \label{Sn-ex}
Let $G=S_n$ for $n \geq 3$, and let $H=S_{n-1}$, the stabilizer of a
point in the usual action of $S_n$ on $n$ points. Let $J$ be the
cyclic group generated by any $n$-cycle in $S_n$. Then $J$ is regular
on the $n$ points, and hence is complementary to $H$. Thus, by Lemma
\ref{untangle}, any Galois extension with group $\Gamma = S_{n-1} \times C_n$  
admits a Hopf Galois structure of type $S_n$. If $n \geq 5$ then $S_n$
has the simple group $A_n$ as a composition factor, whereas $\Gamma$
does not. 
\end{example}

\begin{example} \label{An-ex}
Let $G=A_n$ with $n \geq 5$, odd, and let $H=A_{n-1}$, the stabilizer
of a point. As $n$ is odd, $A_n$ contains an $n$-cycle, and again this
generates a complement $J$ to $H$. Thus any Galois extension with
group $ \Gamma = A_{n-1} \times C_n$ admits a Hopf-Galois structure of
type $A_n$. 
\end{example}

For even $n$, the case $G=A_n$ is a little more involved. 

\begin{lemma} \label{An-gen}
Let $n \geq 4$ with $n \not \equiv 2 \bmod 4$. Write $n=2^e m$ with $m
\geq 1$ odd. Then any Galois extension with group $\Gamma = A_{n-1}
\times C_2^e \times C_m$ admits a Hopf-Galois structure of type
$A_n$.
\end{lemma}
\begin{pf} 
When $n$ is odd (so $e=0$), the result follows from Example
\ref{An-ex}. So suppose $e \geq 2$. We view $A_n$ as permuting the
elements of the finite group $B=C_2^e \times C_m$ of order $n$. We
write $B$ additively. The group $\lambda(B)$ of left translations by
elements of $B$ is a regular subgroup of $\Perm(B) \cong S_n$. 
If $0 \neq v \in C_2^e$, then left translation by $(v,0) \in B$ swaps
the elements of $B$ in pairs. Thus, as a product of disjoint cycles,
$\lambda(v,0)$ consists of $n/2$ transpositions. As $n/2$ is even,
$\lambda(v,0)$ is an even permutation. Similarly, $\lambda(0,1)$ is a
product of $2^e$ $m$-cycles, which is again an even permutation as $m$
is odd. But $B$ is generated by the elements $(v,0)$ for $0 \neq v \in
C_2^e$ and $(0,1)$, so $\lambda(B)$ lies in $A_n$. Thus $J=\lambda(B)$ is 
a complementary subgroup to $A_{n-1}$ in $A_n$, and we may again apply
Lemma \ref{untangle}.
\end{pf}

\begin{remark}
If $n \equiv 2 \bmod 4$, then there is no complementary subgroup $J$
to $A_{n-1}$ in $A_n$. To see this, observe that $J$ would have to be a
regular subgroup, so that any element of $J$ of order $2$ would
consist of $n/2$ transpositions, and therefore could not be in $A_n$.
\end{remark}
  
We now turn to the proof of Theorem \ref{index}. We will need the
following result of Guralnick \cite[Theorem 1]{Gu}, which depends on
the classification of finite simple groups.

\begin{theorem}[Guralnick] \label{guralnick}
Let $G$ be a finite nonabelian simple group with a subgroup $H$ of
prime-power index, $|G:H|=p^a>1$. Then one of the following holds.
\begin{itemize}
\item[(a)] $G=A_n$ with $n=p^a \geq 5$ and $H=A_{n-1}$.
\item[(b)] $G=\PSL_n(q)$ and $H$ is the stabilizer of a point or
  hyperplane in the action of $G$ on the projective space
  $\PP^{n-1}(\F_q)$; here $|G:H|=(q^n-1)/(q-1)=p^a$. 
\item[(c)] $G=\PSL_2(11)$ and $H = A_5$. 
\item[(d)] $G$ is the Matthieu group $M_{23}$ or $M_{11}$, and
  $H=M_{22}$ or $M_{10}$, respectively. 
\item[(e)] $G=\PSU_4(2) \cong \PSp_4(3)$ and $H$ is a maximal
  parabolic subgroup of $\PSU_4(2)$ of index $27$. 
\end{itemize}
Moreover, $H$ is a Hall $p'$-subgroup of $G$ except in case (a) with
$n=p^a>p$ and in case (e).
\end{theorem}

\begin{pf-of}{Theorem \ref{index}}

We consider the various cases in Theorem \ref{guralnick}. In case (a),
the result follows from Lemma \ref{An-gen}. In cases (b), (c) and (d),
the group $H$ is a Hall $p'$-subgroup. Taking $J$ to be a Sylow
$p$-subgroup of $H$, so that $H$, $J$ are complementary subgroups in
$G$, the result follows from Lemma \ref{untangle}. 

It remains to consider case (e), where $H$ is not a Hall
$p'$-subgroup. Let $G=\PSU_4(2)$. As the factor $d$ in Table
\ref{lie-class-autom} is $(n,q+1)=(4,3)=1$ in this case,
$G \cong \SU_4(2)$. Let $V=\F_4^4$, and label the standard ordered basis of
row vectors as $e_1$, $f_1$, $e_2$, $f_2$. We endow $V$ with the
sesquilinear form $( \cdot, \cdot)$ where
$$ (e_i,e_j)=(f_i,f_j)=0, \qquad (e_i,f_j)=\delta_{ij}. $$
Then we may take $G$ to be the subgroup of $\SL_4(4)$ which preserves
this form. More concretely, $G$ consists of the matrices $M$ over
$\F_4$ of determinant $1$ such that $MF\overline{M} = F$, where
$$ F = \left( \begin{array}{cccc} 0 & 1 & 0 & 0 \\ 1 & 0 & 0 & 0 \\ 0
  & 0 & 0 & 1 \\ 0 & 0 & 1 & 0 \end{array} \right) $$ 
and, for any matrix $M$, we write $\overline{M}$ for the transpose of
the matrix obtained by applying the involution $x \mapsto
\overline{x}=x^2$ of $\F_4$ to each entry of $M$. Without loss of
generality, we take $H$ to the stabilizer in $G$ of the maximal
isotropic subspace $W=\F_4 e_1 + \F_4 e_2$ of $V$. Let $X$ be the set
of all $2$-dimensional isotropic subspaces of $V$. Then $G$ acts
transitively on $X$, and a calculation confirms that indeed
$|X|=27$. Thus we need to show that there is a subgroup of $G$ which acts
regularly on $X$.  Such a group will have order $27$, whereas a Sylow
$3$-subgroup of $G$ has order $81$.

Let $\omega \in \F_4$ with $\omega^2+\omega+1=0$, and consider the
matrices 
$$ A= \left( \begin{array}{cccc} 
   1 & \omega & 1 & \omega^2 \\
   \omega^2 & 1 & \omega & 1 \\
   1 & \omega^2 & 1 & \omega \\
   \omega^2 & \omega & \omega & \omega^2 \end{array} \right), \qquad
   B = \left( \begin{array}{cccc}
   1 & 1 & 0 & 0 \cr
   1 & 0 & 0 & 0 \cr
   0 & 0 & 1 & 0 \cr
   0 & 0 & 0 & 1 \end{array} \right) . $$
One verifies by direct calculation that $A$, $B \in G$ and that
$A^9=B^3=I \neq A^3$ and $BA=A^4B$. Thus $A$ and $B$ generate a
subgroup $J$ of $G$ of order $27$. Moreover, the images of $W$ under
the $10$ matrices $A^m$, $0 \leq m \leq 8$ and $B$ are all
distinct. (The image of $W$ in each case is the $\F_4$-span of rows
$1$ and $3$ of the matrix.) Since the size of the orbit of $W$ under $J$ must
be a factor of $|J|$, it follows that $J$ is transitive, and hence
regular, on $X$ as required. 
\end{pf-of}

\begin{pf-of}{Corollary \ref{sol-insol}}

\noindent (i) This follows from Theorem \ref{index} (or Example
\ref{An-ex}) on taking $G=A_5$ and $H=A_4$.

\noindent (ii) Let $G=\PSL_3(2)\cong \PSL_2(7)$ of order $168$. The
stabilizer $H$ of a point (or line) in the projective plane
$\PP^2(\F_2)$ has index $7$ and order $24$. By \cite[8.27]{Hup}, we
have $H \cong S_4$. The result then follows from Lemma \ref{untangle}.

\noindent (iii) Let $p=2^e-1$ be a Mersenne prime, and let 
 $G=\PSL_2(p)$. 
The stabilizer $H$ of a point in $\PP^1(\F_p)$ under
the action of $G$ has index $p+1=2^e$ and order
$\frac{1}{2}p(p-1)$. In particular, $H$ is a Hall $2'$-subgroup of
$G$, so a Sylow $2$-subgroup $J$ will be complementary to $H$.  
From \cite[8.27]{Hup}, $H \cong F_{\frac{1}{2}p(p-1)}$ and $J \cong D_{2^e}$, so
again Lemma \ref{untangle} gives the required result.
\end{pf-of}
 
\bigskip

{\sc \noindent Acknowledgement:} The author would like to thank Griff Elder and
Henri Johnston for helpful comments on an earlier version of this
paper.

\end{document}